\documentclass[a4paper,12pt]{article}
\usepackage{comment}
\usepackage{cite}
\usepackage{amsmath}
\usepackage{amssymb}
\usepackage{amsfonts}
\usepackage[T1]{fontenc}
\usepackage[utf8]{inputenc}
\usepackage{graphicx}
\usepackage{fancyhdr}
\usepackage{float}
\usepackage{xcolor}
\usepackage{authblk}
\usepackage{mathrsfs}
\usepackage{empheq}
\usepackage[hyphens]{url}
\usepackage{hyperref} 
\usepackage[]{breakurl}

\usepackage{geometry}
\DeclareMathOperator{\Tr}{Tr}

\begin{document}

\title{A derivation of Dickson polynomials using the Cayley-Hamilton theorem}

\author[$\dagger$]{Jean-Christophe {\sc Pain}$^{1,2,}$\footnote{jean-christophe.pain@cea.fr}\\
\small
$^1$CEA, DAM, DIF, F-91297 Arpajon, France\\
$^2$Universit\'e Paris-Saclay, CEA, Laboratoire Mati\`ere en Conditions Extr\^emes,\\ 
F-91680 Bruy\`eres-le-Ch\^atel, France
}

\date{}

\maketitle

\begin{abstract}
In this note, the first-order Dickson polynomials are introduced through a particular case of the expression of the trace of the $n^{\mathrm{th}}$ power of a matrix in terms of powers of the trace and determinant of the matrix itself. The technique relies on the Cayley-Hamilton theorem and its application to the derivation of formulas due to Carlitz and to second-order Dickson polynomials is straightforward. Finally, generalization of Dickson polynomials over finite fields and multivariate Dickson polynomials are evoked as potential avenues of investigation in the same framework.
\end{abstract}

\section{Introduction}

The Dickson polynomials are particular cases, for a quadratic equation, of the Waring formulas \cite{Pena_card}, which are expressions for the sums of the $k^{\mathrm{th}}$ power of the roots $x_1, x_2, \cdots, x_n$ of a polynomial equation
\begin{equation}
    x^n+v_1x^{n-1}+\cdots+v_n=0,
\end{equation}
\emph{i.e.},
\begin{equation}
    s_k=\sum_{i=1}^nx_i^k,
\end{equation}
which can be expressed as
\begin{equation}\label{trois}
    s_k=k\sum_{\{r_i\}} (-1)^{\sum_{i=1}^nr_i}\left[\frac{\left(\sum_{i=1}^nr_i-1\right)!}{\prod_{i=1}^nr_i!}\right]\prod_{i=1}^nv_i^{r_i},
\end{equation}
the sum running over the whole set of integers $r_1, r_2, \cdots, r_n \geq 0$ such that $\sum_{i=1}^n i\,r_i = k$ \cite{Dickson1897,Dickson1922}.

In his PhD thesis, Dickson studied a class of polynomials of the form \cite{Pena_card}:
\begin{equation}\label{quatre}
    x^n+n\sum_{i=1}^{\frac{n-1}{2}}\frac{(n-i-1)\cdots (n-2i+1)}{i!}\,b^i\,x^{n-2i},
\end{equation}
over finite fields where $n$ is odd and $b$ a real parameter. Note that expressions (\ref{trois}) and (\ref{quatre}) have coefficients which involve dividing by an integer.  That is not always a meaningful thing to do over a field of characteristic $p>0$, in case the number being divided by is a multiple of $p$. Actually, all the rational numbers occurring in these equations are actually integers, so that the expressions make sense in positive characteristic.

Schur named these polynomials in honor of Dickson and noticed that they are related to the Chebyshev polynomials. Schur defines an integer $n$ to be a ``Dickson number'' if all degree-$n$ polynomials in $\mathbb{Z}[x]$ which permute $\mathbb{Z}/p\mathbb{Z}$ for infinitely many primes $p$ are compositions of linear functions, $x^k$ monomials with $k$ odd, and Dickson polynomials of degree coprime to 6. The main result of Schur's 1923 paper is that all primes $n$ are Dickson numbers.  He claims that he will prove in a subsequent paper that a positive integer $n$ is a Dickson number if it satisfies (C): for each composite divisor $d$ of $n$, every primitive subgroup of the symmetric group $S_d$ which contains a $d$-cycle is doubly transitive. He then notes that Burnside proved that property holds when $n$ is a prime power, so that Schur's announced forthcoming result implies that all prime powers are Dickson numbers. Moreover, Schur writes that he will also prove in a future paper that such a result holds when $n$ is a product of two prime powers, and thus his announced forthcoming result implies that all products of two prime powers are Dickson numbers. Notably, Schur never says anything about whether all positive integers are Dickson numbers. However, a proof of the result announced by Schur, using the arguments from his 1923 paper, was given by Peter M\"uller in 1997 \cite{Muller1997}. Moreover, in a 1933 paper \cite{Schur1933}, Schur proved that every positive integer $n$ satisfies (C).  Note that Schur's 1933 paper contains a remark that in a subsequent paper he will give number-theoretic consequences of his group-theoretic theorem.  

Only polynomials with integer coefficients that induce permutations of the integers {\it modulo} $p$ for infinitely many primes $p$ are compositions of linear polynomials, power polynomials $x^n$, and Dickson polynomials. It is commonly admitted that the polynomials were rediscovered by Brewer in 1960 in his work on the now eponymous sums. They are therefore sometimes (although rarely) referred to as ``Brewer'' polynomials. In mathematics, Brewer sums are finite character sums related to Jacobsthal sums \cite{Brewer1961,Berndt1979}. The Brewer sum is given by
\begin{equation}
    \Lambda_{n}(a)=\sum _{x{\bmod {p}}}\left(\frac{D_{n+1}(x,a)}{p}\right),
\end{equation}
where $D_n$ is the Dickson polynomial \cite{Dickson1897,Lidl1993,Jammes2013} determined by
\begin{equation}\label{recu}
    D_{0}(x,a)=2,\quad D_{1}(x,a)=x,\quad D_{n+1}(x,a)=xD_{n}(x,a)-aD_{n-1}(x,a)
\end{equation}
and $\left({\frac {a}{p}}\right)$ is the usual Legendre symbol\footnote{Let $p$ be an odd prime number. An integer $a$ is a quadratic residue modulo $p$ if it is congruent to a perfect square modulo $p$ and is a quadratic nonresidue modulo $p$ otherwise. The Legendre symbol is a function of $a$ and $p$ defined as
\begin{equation}
    \left({\frac {a}{p}}\right)={\begin{cases}1&{\text{if }}a{\text{ is a quadratic residue modulo }}p{\text{ and }}a\not \equiv 0{\pmod {p}},\\-1&{\text{if }}a{\text{ is a quadratic nonresidue modulo }}p,\\0&{\text{if }}a\equiv 0{\pmod {p}}.\end{cases}}
\end{equation}
}.

Dickson permutation polynomials have been used in cryptography and as a key exchange protocol, in the transmission of information in a secure way \cite{Pena2014}. The Dickson cryptosystem is more general than the RSA (Rivest-Shamir-Adleman) cipher, since for the Dickson scheme the modulus $n$ needs not be squarefree, but can be an arbitrary positive integer with at least two prime factors. It is worth mentioning that Dickson polynomials have also been used for primality testing in number theory \cite{Lidl1993}.

\section{Matrix derivation of the Waring formula and Dickson polynomials}

In this section we will express $(x^n+a^n/x^n)$ in terms of powers of $(x+a/x)$. Let us consider the matrix
\begin{equation}
    M=\left(
    \begin{array}{cc}
    m_{11} & m_{12}\\
    m_{21} & m_{22}\\
    \end{array}
    \right).
\end{equation}
We have
\begin{equation}\label{tra}
    \Tr\left(M^n\right)=\sum_{k=0}^{\lfloor \frac{n}{2}\rfloor}(-1)^k\frac{n}{n-k}\binom{k}{n-k}\left(\Tr M\right)^{n-2k}\left(\det M\right)^k,
\end{equation}
where $\lfloor x \rfloor$ represents the integer part of $x$. Equation (\ref{tra}) can be proven by induction. In the $n=2$ case:
\begin{equation}
    \sum_{k=0}^1(-1)^k\frac{2}{2-k}\binom{k}{2-k}\left(\Tr M\right)^{2-2k}\left(\det M\right)^k=\left(\Tr M\right)^2-2\det M=\Tr\left(M^2\right).
\end{equation}
because of the Cayley-Hamilton theorem \cite{Protat,Householder2006}, which gives
\begin{equation}\label{cy1}
    M^2-(\Tr M)M+(\det M)I_2=0,
\end{equation}
where $I_2$ is the $2\times 2$ identity matrix. Multiplying the latter expression by $M$ gives
\begin{equation}
    M^3-(\Tr M)M^2+(\det M)M=0,
\end{equation}
and thus
\begin{equation}\label{cy2}
    \Tr\left(M^3\right)=(\Tr M)(\Tr M^2)-(\det M)(\Tr M)
\end{equation}
or, inserting Eq. (\ref{cy1}) in Eq. (\ref{cy2}):
\begin{equation}
    \Tr\left(M^3\right)=(\Tr M)\left[(\Tr M)^2-2(\det M)\right]-(\det M)(\Tr M)=(\Tr M)^3-3(\det M)(\Tr M)
\end{equation}
and thus, for $n=3$:
\begin{equation}
    \sum_{k=0}^1(-1)^k\frac{3}{3-k}\binom{k}{3-k}\left(\Tr M\right)^{3-2k}\left(\det M\right)^k=\left(\Tr M\right)^3-3(\Tr M)(\det M)=\Tr\left(M^3\right).
\end{equation}
Let us assume that the formula is true at ranks $n$ and $n+1$, for $n\geq 2$. We have
\begin{equation}
    M^{n+2}=(\Tr M)M^{n+1}-(\det M)M^n=0,
\end{equation}
and 
\begin{align}
    \Tr\left(M^{n+2}\right)&=(\Tr M)\sum_{k=0}^{\lfloor \frac{n+1}{2}\rfloor}(-1)^k\frac{n+1}{n+1-k}\binom{k}{n+1-k}\left(\Tr M\right)^{n+1-2k}\left(\det M\right)^k\nonumber\\
    & -(\det M)\sum_{k=0}^{\lfloor \frac{n}{2}\rfloor}(-1)^k\frac{n}{n-k}\binom{k}{n-k}\left(\Tr M\right)^{n-2k}\left(\det M\right)^k
\end{align}
or equivalently
\begin{align}
    \Tr\left(M^{n+2}\right)&=\sum_{k=0}^{\lfloor \frac{n+1}{2}\rfloor}(-1)^k\frac{n+1}{n+1-k}\binom{k}{n+1-k}\left(\Tr M\right)^{n+2-2k}\left(\det M\right)^k\nonumber\\
    & -\sum_{k=0}^{\lfloor \frac{n}{2}\rfloor}(-1)^k\frac{n}{n-k}\binom{k}{n-k}\left(\Tr M\right)^{n-2k}\left(\det M\right)^{k+1}.
\end{align}
If $n=2p$, the latter equation can be recast into
\begin{align}
    \Tr\left(M^{n+2}\right)&=\sum_{k=0}^{p}(-1)^k\frac{n+1}{n+1-k}\binom{k}{n+1-k}\left(\Tr M\right)^{n+2-2k}\left(\det M\right)^k\nonumber\\
    & -\sum_{k=0}^{p}(-1)^k\frac{n}{n-k}\binom{k}{n-k}\left(\Tr M\right)^{n-2k}\left(\det M\right)^{k+1}\nonumber\\
    &=(\Tr M)^{n+2}+\sum_{k=1}^{p+1}(-1)^k\left[\frac{n+1}{n+1-k}\binom{k}{n+1-k}\right.\nonumber\\
    &\left.+\frac{n}{n+1-k}\binom{k-1}{n+1-k}\right]\left(\Tr M\right)^{n+2-2k}\left(\det M\right)^k
\end{align}
and similarly if $n=2p+1$:
\begin{align}
    \Tr\left(M^{n+2}\right)&=\sum_{k=0}^{p+1}(-1)^k\frac{n+1}{n+1-k}\binom{k}{n+1-k}\left(\Tr M\right)^{n+2-2k}\left(\det M\right)^k\nonumber\\
    & -\sum_{k=0}^{p}(-1)^k\frac{n}{n-k}\binom{k}{n-k}\left(\Tr M\right)^{n-2k}\left(\det M\right)^{k+1}\nonumber\\
    &=(\Tr M)^{n+2}+\sum_{k=1}^{p+1}(-1)^k\left[\frac{n+1}{n+1-k}\binom{k}{n+1-k}\right.\nonumber\\
    &\left.+\frac{n}{n+1-k}\binom{k-1}{n+1-k}\right]\left(\Tr M\right)^{n+2-2k}\left(\det M\right)^k
\end{align}
and thus
\begin{equation}
    \Tr\left(M^{n+2}\right)=\Tr\left(M\right)^{n+2}+\sum_{k=1}^{\lfloor\frac{n+2}{2}\rfloor}(-1)^k\left[\frac{(n+1)\binom{k}{n+1-k}+n\binom{k-1}{n+1-k}}{n+1-k}\right]\left(\Tr M\right)^{n+2-2k}\left(\det M\right)^k.
\end{equation}
We have
\begin{align}
    \frac{(n+1)\binom{k}{n+1-k}+n\binom{k-1}{n+1-k}}{n+1-k}&=\binom{k}{n+1-k}+\frac{k}{n-k+1}\binom{k}{n+1-k}\nonumber\\
    &+\binom{k-1}{n+1-k}+\frac{k-1}{n-k+1}\binom{k-1}{n+1-k}\nonumber\\
    &=\binom{k}{n+1-k}+\binom{k-1}{n-k}+\binom{k-1}{n+1-k}+\binom{k-2}{n-k}\nonumber\\
    &=\binom{k}{n+1-k}+\binom{k-1}{n+1-k}+\binom{k-1}{n-k}+\binom{k-2}{n-k}\nonumber\\
    &=\binom{k}{n+2-k}+\binom{k-1}{n+1-k}\nonumber\\
    &=\left(1+\frac{k}{n+2-k}\right)\binom{k}{n+2-k}=\frac{n+2}{n-k+2}\binom{k}{n+2-k}
\end{align}
yielding
\begin{equation}
    \Tr\left(M^{n+2}\right)=\sum_{k=0}^{\lfloor\frac{n+2}{2}\rfloor}(-1)^k\frac{(n+2)}{n-k+2}\binom{k}{n+2-k}\left(\Tr M\right)^{n+2-2k}\left(\det M\right)^k,
\end{equation}
which completes the proof. Setting $m_{11}=x$, $m_{12}=m_{21}=0$ and $m_{22}=y$, we get
\begin{equation}
    M=\left(
    \begin{array}{cc}
    x & 0\\
    0 & y\\
    \end{array}
    \right),
\end{equation}
for which $\Tr\left(M^n\right)=x^n+y^n$ and $\det M=xy$. This yields
\begin{equation}
    x^n+y^n=\sum_{k=0}^{\lfloor\frac{n}{2}\rfloor}(-1)^k\frac{n}{n-k}\binom{k}{n-k}\left(x+y\right)^{n-2k}(xy)^k
\end{equation}
which is sometimes referred to as the Waring formula \cite{Gould1999} or also, for $y=1/x$:
\begin{equation}
    x^n+\frac{1}{x^n}=\sum_{k=0}^{\lfloor\frac{n}{2}\rfloor}(-1)^k\frac{n}{n-k}\binom{k}{n-k}\left(x+\frac{1}{x}\right)^{n-2k}
\end{equation}
or again, more generally, for $y=a/x$:

\begin{equation}
    x^n+\frac{a^n}{x^n}=\sum_{k=0}^{\lfloor\frac{n}{2}\rfloor}(-a)^k\frac{n}{n-k}\binom{k}{n-k}\left(x+\frac{a}{x}\right)^{n-2k},
\end{equation}
where we find the Dickson polynomials of the first kind $(D_{n})_{n\in\mathbb{N}}$, defined as
\begin{equation}
    D_{n}(x,a)=\sum \limits _{k=0}^{\lfloor n/2\rfloor}{\frac{n}{n-k}}{\binom{n-k}{k}}\,(-a)^{k}\,x^{n-2k},
\end{equation}
and satisfying the recurrence relation (\ref{recu}):
\begin{equation}
    D_{n+2}(x,a)=x\,D_{n+1}(x,a)-a\,D_{n}(x,a).
\end{equation}
The first Dickson polynomials are:
\begin{equation}
    \begin{array}{ll}
    D_{0}(x,a)=2\\
    D_{1}(x,a)=x\\
    D_{2}(x,a)=x^{2}-2a\\
    D_{3}(x,a)=x^{3}-3ax\\
    D_{4}(x,a)=x^{4}-4ax^{2}+2a^{2}\\
    D_{5}(x,a)=x^{5}-5ax^{3}+5a^{2}x.
    \end{array}
\end{equation}
One thus has, $\forall(x,a)\in \mathbb{C}^{*}\times\mathbb {C}$ and $n\in \mathbb{N}$:
\begin{equation}\label{solu}
    D_{n}\left(x+{\frac {a}{x}},a\right)=x^{n}+{\frac{a^{n}}{x^{n}}}.
\end{equation}
The Dickson polynomials satisfy, for $(m,n)\in \mathbb {N}^{2}$ : $D_{mn}(x,a)=D_{m}\left(D_{n}(x,a),a^{n}\right)$ and are solutions of the differential equation: $(x^{2}-4a)y''+xy'-n^{2}y=0$. Their generating series reads \cite{Wilf1990}:
\begin{equation}
    \sum_{n\in\mathbb{N}}D_{n}(x,a)\,z^{n}=\frac{2-xz}{1-x\,z+a\,z^{2}}.
\end{equation}
The Dickson polynomials are related to the Chebyshev polynomials of the first kind $(T_{n})_{n\in \mathbb{N}}$ by the relation: $\forall (x,a)\in \mathbb{C}^{2}$, $D_{n}(2ax,a^{2})=2a^{n}\,T_{n}(x)$ (see Appendix \ref{appA} for a proof of Eq. (\ref{solu}) relying on that correspondence). On the other hand, the Dickson polynomials for which $a=0$ are monomials $D_{n}(x,0)=x^{n}$ ; the ones with $a=1$ and $a=-1$ are related to Fibonacci and Lucas polynomials respectively.

\section{Higher-order formulas}

The same procedure can be applied to more than two variables, in order to derive other identities, such as the formula due to Carlitz \cite{Carlitz1978,Carlitz1978b}:
\begin{equation}
    \sum_{i+2j+3k=n}(-1)^j\frac{n}{i+j+k}\frac{(i+j+k)!}{i!j!k!}(x+y+z)^{n-3k}(xy+yz+zx)^k=x^n+y^n+z^n
\end{equation}
summed over all $0\leq i,j,k\leq n$, $n>0$ and where $xyz=1$. One has also
\begin{equation}
    x^n+y^n+z^n=\sum_{k=0}^{\lfloor\frac{n}{3}\rfloor}\frac{n}{n-2k}\binom{k}{n-2k}\left(x+y+z\right)^{n-2k}(xyz)^k,
\end{equation}
where $xy+yz+zx=0$.

\section{Conclusion}

We proposed a simple derivation of the first-order Dickson polynomials, as a particular case of the expression of the trace of the $n^{\mathrm{th}}$ power of a $2\times 2$ matrix in terms of powers of the trace and determinant of the matrix itself. The same idea can be applied in order to obtain the Carlitz expression (thus involving Cayley-Hamilton theorem for $3\times 3$ matrices), and generalized to higher powers. The methodology applies to the Dickson polynomials of the second kind, and possible extensions to generalized and multivariate Dickson polynomials (see Appendix \ref{appA}) would be worth investigating.

\appendix

\section{Recovering $x^n+a/x^n$ using the Chebyshev polynomials as particular cases of Dickson ones}\label{appA}

Let us set $a=b^2$, $x=by$ and $y=e^{i\theta}$. One has
\begin{equation}\label{eq1}
    x^n+\frac{a^n}{x^n}=b^n\,y^n+\frac{b^n}{y^n}=b^n\,(e^{in\theta}+e^{-in\theta})=2b^n\,\cos(n\theta).
\end{equation}
On the other hand,
\begin{equation}\label{eq2}
    D_n\left(x+\frac{a}{x},a\right)=D_n(2by,b^2)=2b^n\,T_n(\cos\theta)=2b^n\cos(n\theta).
\end{equation}
From Eqs. (\ref{eq1}) and (\ref{eq2}), one gets
\begin{equation}
    D_n\left(x+\frac{a}{x},a\right)=x^n+\frac{a^n}{x^n}.
\end{equation}

\section{Generalized Dickson polynomials}\label{appB}

Dickson polynomials over finite fields can be thought of as initial members of a sequence of generalized Dickson polynomials referred to as Dickson polynomials of the $(k + 1)^{\mathrm{th}}$ kind \cite{Filipponi1997,Young2002,Wang2012}. Specifically, for $a\ne 0 \in \mathbb{F}_q$, $\mathbb{F}_q$ being a finite field with $q = p^m$ for some prime $p$ and any integers $n \geq 0$ and $0 \leq k < p$, the $n^{\mathrm{th}}$ Dickson polynomial of the $(k + 1)^{\mathrm{th}}$ kind over $\mathbb{F}_q$, denoted by $D_{n,k}(x,\alpha)$, is defined by $D_{0,k}(x,\alpha)=2-k$ and
\begin{equation}
D_{n,k}(x,a)=\sum_{i=0}^{\left\lfloor{\frac{n}{2}}\right\rfloor}{\frac{n-ki}{n-i}}{\binom{n-i}{i}}\,(-a)^{i}\,x^{n-2i}.
\end{equation}
Note that one has as a special case $D_{n,0}(x,a)=D_n(x,a)$.

The main properties of the Dickson polynomials also generalize. As concerns the recurrence relation, for instance, one has, for $n \geq 2$:
\begin{equation}
D_{n,k}(x,a)=x\,D_{n-1,k}(x,a)-a\,D_{n-2,k}(x,a),
\end{equation}
with the initial conditions $D_{0,k}(x,a)=2-k$ and $D_{1,k}(x,a)=x$. The functional equation for $D_{n,k}$ polynomials reads
\begin{align}\label{solu40}
D_{n,k}\left(y+\frac{a}{y}\right)&={\frac{y^{2n}+k\,a\,y^{2n-2}+\cdots+ka^{n-1}\,y^{2}+a^{n}}{y^{n}}}\nonumber\\
&={\frac{y^{2n}+{a}^{n}}{y^{n}}}+\left({\frac{ka}{y^{n}}}\right){\frac{y^{2n}-{a}^{n-1}\,y^{2}}{y^{2}-a}},
\end{align}
where $y \ne 0$, $y^2 \ne a$. Their generating function is
\begin{equation}
\sum_{n=0}^{\infty}D_{n,k}(x,a)\,z^{n}=\frac{2-k+(k-1)\,x\,z}{1-x\,z+a\,z^{2}}.
\end{equation}
Finally, the Dickson polynomials can be extended to the multivariate case (see Ref. \cite{Mullen2013}, p. 274 and Ref. \cite{Zriaa2023}):
\begin{equation}
    D_n^{(i)}(x_1,\cdots, x_t,a)=s_i(u_1^n, \cdots, u_{t+1}^n)
\end{equation}
for $1\leq i\leq t$, where $x_i=s_i(u_1, \cdots, u_{t+1})$ are elementary symmetric functions and $u_1\cdots u_{t+1}=a$. It satisfies the generating function
\begin{equation}
    \sum_{n=0}^{\infty}D_n(x_1,\cdots,x_t,a)\,z^n=\frac{\sum_{i=0}^t(t+1-i)\,(-1)^i\,x_i\,z^i}{\sum_{i=0}^{t+1}(-1)^i\,x_i\,z^i}
\end{equation}
for $n\geq 0$. The recurrence relation is
\begin{equation}
    D_{n+t+1}^{(1)}-x_1D_{n+t}^{(1)}+\cdots+(-1)^t\,x_k\,D_{n+1}^{(1)}+(-1)^{t+1}\,a\,D_n^{(1)}=0,
\end{equation}
where the initial values are given by $D_0^{(1)}=t+1$ and
\begin{equation}
    D_j^{(1)}=\sum_{r=1}^j(-1)^{r-1}\,x_r\,D_{j-r}^{(1)}+(-1)^j\,(t+1-j)\,x_j,
\end{equation}
for $0<j\leq t$.

At this stage we have to make some comments. The most important property of Dickson polynomials is the functional equation (\ref{solu}), because it can be rewritten as  $D_n(x,a) \circ (x+a/x) = (x+a^n/x) \circ x^n$,  which in some sense says that the Dickson polynomial is a ``coordinate projection'' of $x^n$.  This perspective becomes especially clear when written for Chebyshev polynomials, where it reads $T_n(z) \circ [(z+1/z)/2] = [(z+1/z)/2] \circ z^n$.  When $z=e^{i \theta}$, we have $z = \cos \theta + i \sin \theta$ and $z^n = e^{i n \theta} = \cos(n \theta) + i \sin (n \theta)$, and the equation becomes the familiar identity  $T_n(\cos \theta) = \cos (n \theta)$.  Note that in the complex plane, $\cos \theta$ is the $x$-coordinate of $z$, and $\cos(n \theta)$ the $x$-coordinate of $z^n$. Thus, $T_n$ is the map on $x$-coordinates corresponding to $z\rightarrow z^n$. Likewise, $D_n(x,a)$ comes from $x^n$ by taking a degree-2 coordinate projection in both the domain and target spaces, via $x+a/x$ and $x+a^n/x$, respectively.  The fact that $D_n(x,a)$ and $T_n(x)$ are coordinate projections of $x^n$ explains all of their interesting properties and applications. There are very few results in the literature about Dickson polynomials of the second and higher kinds. It is true that Eq. (\ref{solu40}) is in some sense a generalization of Eq. (\ref{solu}), but this generalization is not a functional equation, because the right-hand side is not a function of $y^n$.  This is why the higher-kind Dickson polynomials may not have really interesting properties. 

\section*{Acknowledgments}

I would like to thank Michael Zieve for useful comments, especially as concerns Schur's 1923 article as well as the generalized Dickson polynomials.

\end{document}